
\input amstex
\magnification=\magstep1
\documentstyle{amsppt}

\redefine\sp{\operatorname{\;sp}}
\define\sa{\operatorname{sa}}
\define\Tr{\operatorname{\;Tr}}
\define\ux{\underline{x}}

\define\ul{\underline{\lambda}}
\define\uI{\underline{I}}

\define\tA{\widetilde {A}}
\define\tC{\widetilde {C}}
\define\tT{\widetilde {T}}

\overfullrule=0pt
\TagsOnRight

\topmatter
\title
Convex multivariable trace functions
\endtitle
\rightheadtext{Convex multivariable trace functions }
\author
Elliott H. Lieb \& Gert K. Pedersen
\endauthor
\date{October 10, 2001}
\enddate
\address{Departments of Mathematics and Physics, Princeton University, 
P.O. Box 708,\linebreak
Princeton NJ 08544-0708, USA 
\& Department of Mathematics,  University of Copenhagen, 
Universitetsparken 5, DK-2100 Copenhagen \O, Denmark} 
\endaddress
\email{lieb\@math.princeton.edu \& gkped\@math.ku.dk}\endemail

\abstract{For any densely defined, lower semi-continuous trace $\tau$ 
on a $C^*-$algebra $A$ with mutually commuting $C^*-$subalgebras 
$A_1,A_2, \dots A_n$, and a convex  function $f$ of $n$ variables, we
give a short proof of the fact that the function 
$(x_1, x_2, \dots , x_n) @>>> \tau (f(x_1, x_2, \dots , x_n))$ is 
convex on the space $\bigoplus_{i=1}^{\;\; n} (A_i)_{\sa}$. If
furthermore the function $f$ is log-convex or root-convex, so is
the corresponding trace function. We also introduce a generalization of 
log-convexity and root-convexity called $\ell-$convexity, show how it
applies to traces, and give some examples. In particular we show that the
Kadison-Fuglede determinant is concave and that the trace of an operator mean
is always dominated by the corresponding mean of the trace values.}   
\endabstract


\endtopmatter
 
\document

\footnote""{\copyright 2001 by the authors. This paper may be reproduced,
in its entirety, for non-commercial purposes.}

\subhead{1. Introduction}\endsubhead 
The fact that several important concepts in operator theory, in quantum
statistical mechanics (the entropy, the relative entropy, Gibbs free
energy), in engineering and in economics involve the trace of a function
of a self-adjoint operator has motivated a considerable amount of abstract
research about such functions in the last half century. An
important subset of questions involve the convexity of
trace functions with respect to their argument. 

The convexity of the function
$x@>>> \Tr (f(x))$, when $f$ is a convex function of one variable and
$x$ is a self-adjoint operator, was known to von Neumann, cf\.
\cite{{\bf 12}, V.3. p. 390}. An early proof for $f(x) = \exp(x)$ can be
found, e\.g\. in \cite{{\bf 19}, 2.5.2}. A proof given by the first author 
some time ago describes the trace $\Tr (f(x))$, where $f$ is convex, as a 
supremum over all possible choices of orthonormal bases of the Hilbert space
of the sum of the values of $f$ at the diagonal elements of the matrix for $x$.
This proof was communicated to B. Simon, who used the method  to give an 
alternative proof of the second Berezin-Lieb inequality in \cite{{\bf 21},
Theorem 2.4}, see also \cite{{\bf 22}, Lemma II.10.4}. Simon only considers
the exponential function, but the argument is valid for any convex function,
cf\. \cite{{\bf 10}, Proposition 3.1}. The general case for an arbitrary trace
on a von Neumann algebra was established by D. Petz in \cite{{\bf 17}, Theorem
4} using the theory of spectral dominance (spectral scale).

The basic fact, for one variable $x$ and a positive convex function $f$, is
that
$$
 \sum\nolimits_j f( (\phi_j, x \ \phi_j)) \le \Tr (f(x))  \,,    
\tag{1}
$$
where the sum -- finite or not -- is over any orthonormal basis. Equality is
obviously achieved if the basis is the set of eigenvectors of $x$. Thus,
$$ \Tr (f(x)) = \sup_{\{\phi \}} \sum\nolimits_j f( (\phi_j, x \ \phi_j)) \ .
\tag2
$$
Essentially, the proof is the following: If $\{\phi_j\}$ is the
eigenvector basis (with eigenvalues $\lambda_j$) and $\psi_j$ is some
other basis,  then $\psi_j = \sum_k C_{jk} \phi_k$, and the
coefficients of the unitary matrix $C$ satisfy $\sum_j  |C_{jk}|^2 =
1 = \sum_k  |C_{jk}|^2$. Then $(\psi_j, x \ \psi_j) = \sum_k |C_{jk}|^2
\lambda_k $ and, by Jensen's inequality, $f(\sum_k  |C_{jk}|^2 \lambda_k)
\leq \sum_k |C_{jk}|^2 f(\lambda_k)$. Now, summing on $j$ we obtain (1).

Equation (2) implies the convexity of $x@>>> \Tr (f(x))$, because any
supremum of convex functions is convex. Moreover, if $f(x) = \exp (x)$
then one sees that $\Tr(\exp(x)) $ is log-convex 
(i.e. $\log (\Tr (\exp(x)))$ is convex), because an ordinary sum of the type
$\sum \exp(a_j)$, with $a_j$ in $\Bbb R$, is log-convex.
Similarly, if $f(x) = |x|^p$ (with $p\geq 1$) we see that $x@>>> \left(
\Tr (|x|^p) \right)^{1/p}$ is convex. In particular, the Schatten $p-$norms are
subadditive.

Even more is true. If $f(x)= \exp (g(x))$ and $g$ is convex, then
$x@>>> \Tr (f(x))$ is log-convex. Similarly, if $f(x)= |g(x)|^p$
then $x@>>> \left( \Tr (f(x)) \right)^{1/p}$ is convex.

A natural question that arises at this point is this: Are there other
pairs of functions  $e,\,\ell$ of one real variable, beside the
pairs $\exp, \,\log$ and $|t|^p,\, |t|^{1/p}$, for which 
$x@>>> \ell \left( \Tr (f(x))\right) $ is convex
whenever $f$ is $\ell-$convex, i.e. $f(x)= e\left(g(x)\right)$ and
$g$ is convex?  In the second part of our paper we answer this question
completely and give a few examples, which we believe to have some potential 
value.  But first we turn to the question of generalizing (2) to functions of
several variables.

We start with a function $f(\ul)$ of $n$ real variables (with $\ul =
(\lambda_1, \lambda_2, \dots , \lambda_n)$). Next we replace the real
variables $\lambda_j$, by operators $x_j$, similar
to the one-variable case. An immediate problem that arises is how to
define  $f(\ux)$ in this case. The spectral theorem, which was used
in the one-variable case, fails here unless the $x_j$'s commute with
each other. Therefore, we restrict the $x_1, x_2, \dots , x_n$ to lie
in commuting subalgebras  $A_1, A_2, \dots , A_n$, and then $f(\ux)$
and $\Tr (f(\ux))$ are well defined and it makes sense to discuss the
joint convexity of this trace function under the condition that $f$
is a jointly convex function of its arguments. (We do not investigate
the question whether $f(\ux)$ is operator convex -- only the convexity
under the trace.)

More generally, we replace the trace $\Tr$ in a Hilbert space
setting by $\tau$, a densely defined, lower semi-continuous  trace  on
a $C^*-$algebra $A$; i\.e\. a functional defined on the set $A_+$ of
positive elements with values in $[0, \infty]$, such that $\tau (x^*x) 
=\tau (xx^*)$ for all $x$ in $A$. We further assume that $A$ comes equipped
with mutually commuting $C^*-$subalgebras  $A_1,A_2, \dots , A_n$.

The convexity of the function $\ux @>>> \tau (f(\ux))$ on the space of
$n-$tuples in \break $\bigoplus_{i=1}^{\;\; n} (A_i)_{\sa}$ was proved
by F. Hansen for matrix algebras in \cite{{\bf 5}}. (Here $A_{\sa}$
denotes the self-adjoint elements in  $A$.) His result was extended to
general operator algebras by the second author in \cite{{\bf 16}}. Both
these proofs rely on Fr\'echet differentiability and some rather intricate
manipulations with first and second order differentials.

We realized that the argument in (2) will work in this multi-variable
case as well, thereby providing a quick proof of the convexity. The key
observation is that the mutual commutativity of $x_1, x_2, \dots , x_n$
implies that there is {\it one} orthonormal basis  (in the case
that $A$ is a matrix algebra) that simultaneously makes
{\it all} the $x_j$ diagonal. 

For a general $C^*-$algebra (e.g. the algebra of continuous functions
on an interval) we have to find something to take the place of an
orthonormal basis. This something is just the commutative
$C^*-$subalgebra generated by the $n$ commuting elements $x_1, x_2, \dots ,
x_n$. It depends on $\ux$, of course, but that fact is immaterial for
computing the trace for a given $\ux$. 

The  main result in the first part of this paper -- proved in section 6 -- is

\proclaim{2. Theorem  (Multivariable Convex Trace Functions)}
Let f be a continuous convex function defined  on a 
cube $\uI = I_1\times \cdots \times I_n$ in $\Bbb R^n$. If $A_1, \dots
, A_n$ are  mutually commuting  $C^*-$subalgebras of a $C^*-$algebra
$A$ and  $\tau$ is a finite  trace on $A$, then the function
$$
(x_1, \dots , x_n) @>>> \tau (f(x_1, \dots , x_n))\,, 
\tag3
$$
defined on commuting $n-$tuples such that $x_i\in (A_i)_{sa}^{I_i}$
for each $i$,  is convex on $\bigoplus (A_i)_{\sa}^{I_i}$. (Here, 
$(A_i)_{\sa}^{I_i}$ denotes the set of self-adjoint elements in
$A_i$ whose spectra are contained in $I_i$.)

If $\tau$ is only densely defined, but lower semi-continuous, the
result still holds if $f\ge 0$, even though the
function may now attain infinite values.
\endproclaim

\medskip

In the second part we explore the natural generalization of the concept
of log\-convexity mentioned before and explained in detail in section
7. We find a necessary and sufficient condition on a {\it concave}
function $\ell$ that ensures that $\ell-$convexity of a function $f =
e\circ g$, with $e=\ell^{-1}$ and $g$ convex, implies $\ell-$convexity
of the function $\ux @>>> \tau (f(\ux))$ for a tracial state $\tau$ on
a $C^*-$algebra $A$. (A {\it tracial  state} is a trace satisfying $\tau
(\bold 1) = 1$.) The main result there -- proved in section 9 -- is

\proclaim{3. Theorem ($\ell-$Convex Trace Functions)} Let f be a
continuous function defined  on a cube $\uI = I_1\times \cdots \times
I_n$ in $\Bbb R^n$, and assume furthermore that $f$ is $\ell-$convex
relative to a pair of functions  $e, \ell$ as described in section 7, where
$\ell'/\ell''$ is convex. If  $A_1, \dots , A_n$ are  mutually
commuting  $C^*-$subalgebras of a $C^*-$algebra $A$ and  $\tau$ is a
tracial state on $A$, then the function
$$
(x_1, \dots , x_n) @>>> \tau (f(x_1, \dots , x_n))\,, 
\tag4
$$
defined on commuting $n-$tuples such that $x_i\in (A_i)_{sa}^{I_i}$
for each $i$, is also $\ell-$convex on $\bigoplus
(A_i)_{\sa}^{I_i}$. If moreover $\ell'/\ell''$ is homogeneous and
$f\ge 0$ the result holds for any densely defined, lower
semi-continuous trace $\tau$ on $A$. 
\endproclaim

\medskip

In sections 4, 5 and 7 we set up some necessary machinery, whereas the key 
lemmas are in sections 6 and 8. 

The third part of the paper, sections 9--22, consists of examples where we
apply the preceeding results. In particular we investigate the n-fold harmonic
mean of positive operators and show how its trace behaves under certain
concave transformations. As a corollary we prove in Proposition 23
that the trace of any mean (in the sense of Kubo and Ando \cite{{\bf
8}}) is dominated by the corresponding mean of the trace values.

Throughout the paper we have chosen a $C^*-$algebraic setting with densely
defined, lower semi-continuous traces, this being the more general theory. We
might as well have developed the theory for von Neumann algebras with normal,
semi-finite traces; in fact we need this more special setting in the proof of
Lemma 6. However, the Gelfand-Naimark-Segal construction effortlessly
transforms the $C^*-$algebra version into the von Neumann algebra setting, so
there is no real difference between the two approaches.
  
\bigskip 

\subhead{4. Spectral Theory}\endsubhead We consider a $C^*-$algebra
$A$ of operators on some Hilbert space $\frak H$ and mutually commuting 
$C^*-$subalgebras $A_1,\cdots,A_n$, i\.e\. $A_i \subset A_j$' for all
$i\ne j$. For each interval $I_i$  we let $(A_i)_{\sa}^{I_i}$ denote
the convex set of self-adjoint elements in 
$A_i$ with spectra contained in $I_i$. If   $\uI=I_1\times\cdots\times
I_n\subset\Bbb R^n$ and $f$ is a continuous function on $\uI$ we
can for each $\ux=(x_1,\cdots,x_n)$ in $\bigoplus (A_i)_{\sa}^{I_i}$
define an element $f(\ux)$ in $A$. To see this, let $x_i=\int \lambda
dE_i(\lambda)$ be the spectral resolution of  $x_i$ for $1\le i \le
n$.  Since the $x_i$'s commute, so do their  spectral measures. We can
therefore define the product spectral measure $E$ on $\uI$ by
$E(S_1\times\cdots\times S_n)=E_1(S_1)\cdots E_n(S_n)$ and then write   
$$ 
f(\ux)=\int f(\lambda_1,\cdots,\lambda_n)\,dE
(\lambda_1,\cdots,\lambda_n)\,. \tag5
$$
Of course, if $f$ is a polynomial in the variables 
$\lambda_1, \dots ,\lambda_n$ we simply find $f(\ux)$ by replacing each 
$\lambda_i$ with $x_i$. The map 
$f\to f(\ux)$ so obtained is a $^*-$homomorphism of
$C(\uI)$ into $A$ and generalizes the ordinary spectral mapping
theory for a single (self-adjoint) operator. In particular, the
support of the map (the smallest closed set $S$ such that $f(\ux)=0$
for any function $f$ that vanishes off $S$) may be regarded as the
``joint spectrum'' of the elements $x_1,\dots , x_n$.

This theory applies readily in the situation where
$A=A_1\otimes\cdots\otimes A_n$, but, curiously enough, the tensor product
structure (used extensively in \cite{{\bf 5}} and \cite{{\bf 16}}) is not
needed in our arguments.

\bigskip 

\subhead{5. Conditional Expectations}\endsubhead Let $\tau$ be a fixed,
densely defined, lower semi-continuous trace on $A$, and let $C$ be a fixed
commutative $C^*-$subalgebra of $A$. By Gelfand theory
we know that each commutative  $C^*-$sub\-algebra of $A$ has the
form $C_0(T)$ for some locally  compact  Hausdorff space $T$. Note now
that if $y\in C_+$, the positive part of $C$, and has compact support as
a function on $T$, then
$y=yz$ for some $z$ in $C_+$. Since the minimal dense ideal $K(A)$ of
$A$ is generated (as a hereditary $^*-$subalgebra) by elements $a$ in
$A_+$ such that $a=ab$ for some $b$ in $A_+$, cf\. \cite{{\bf 14},
  5.6.1}, and since $\tau$ is densely defined, hence finite on $K(A)$,
it follows that $\tau (y) < \infty$. Restricting $\tau$ to $C$ we
therefore obtain a unique Radon measure $\mu_C$ on $T$ such that   
$$
\int y(t)\,d\mu_C (t) = \tau (y)\,, \qquad y\in C\,, \tag6
$$
cf\. \cite{{\bf 15}, Chapter 6}. Furthermore, if $x\in A_+$ the
positive functional $y \to \tau (yx)$ on $C$  determines a unique
Radon measure on $T$ (by the Riesz representation theorem)  which is
absolutely continuous with respect to \, $\mu_C$, in fact dominated by a 
multiple  of $\mu_C$ (by the Cauchy-Schwarz inequality). By the Radon-Nikodym
theorem, there is a positive function $\Phi (x)$ in $L^\infty_{\mu_C} (T)$ such
that  
$$
\int y(t) \Phi (x)(t)\,d\mu_C (t) = \tau (yx)\,, \qquad y\in C\,. 
\tag7
$$
Extending by linearity, this defines a map $\Phi$ from $A$ to 
$L^\infty_{\mu_C} (T)$ which is linear, positive, norm decreasing  (and
unital if both $A$ and $C$ have the same unit).
Moreover, $\Phi (y) =y$ almost everywhere if $y \in C$ (with respect
to the natural  homomorphism of $C = C_0(T)$ into $L^\infty_{\mu_C} (T)$ ).

When $\tau$ is faithful and $C$ and $A$ are von Neumann algebras, the
map $\Phi$ is a classical example of a conditional expectation,
cf\. \cite{{\bf 7}, Exercise 8.7.28}.

\bigskip 

\proclaim{6. Lemma} With notations as in sections 4 and  5, take 
$\ux = (x_1, \dots , x_n)$ in $\bigoplus_{i=1}^{\;\;n} (A_i)_{\sa}^{I_i}$ 
and let $f$ be a continuous, real and convex function on 
$\uI=I_1\times \cdots I_n$. If $\tau$ is unbounded we assume,
moreover,  that $f$ is positive. For each commutative $C^*-$subalgebra $C$ of
$A$ set
$$
\varphi_C (f, \ux) =\int f \bigl(\Phi (x_1)(t), \dots ,
\Phi(x_n)(t)\bigr)\,d\mu_C (t) \,.  \tag8
$$
Then $\varphi_C (f, \ux)\le \tau (f(\ux))$, with equality whenever $x_i \in C$
for all $i$. 
\endproclaim 

\demo{Proof} If $x_i \in C$ for all $i$ then also $f(\ux)\in C$, and
since $\Phi (x_i)=x_i$ almost everywhere we get by (6) that
$$
\varphi_C (f, \ux) = \int f\bigl(x_1(t), \dots , x_n(t)\bigr)\,d\mu_C (t) 
= \int f(\ux)(t)\,d\mu_C (t) = \tau (f(\ux))\,.  \tag9
$$
 
To prove the inequality for a general $n-$tuple we first assume that
$\tau$ is a normal, semi-finite trace on a von Neumann algebra $M$, 
\cite{{\bf 7}, 8.5} or \cite{{\bf 20}, 2.5.1}, and put
$$
M_\tau = \{ x\in M\mid \tau (|x|) < \infty\}\,.
\tag10
$$
We then assume that $A=M_\tau^=$, the norm closure  of $M_\tau$, so
that $A$ is a two-sided ideal in $M$. We further assume that each
$A_i$ is relatively weakly closed in $A$, 
i\.e\. $A_i=A\cap M_i$ for some von Neumann subalgebra $M_i$ of
$M$. This has the effect that every self-adjoint element $x_i$ in
$A_i$ can be approximated in norm by an element $y_i$ in $A_i$ with
finite spectrum and $\sp (y_i)  \subset \sp (x_i)$. Under these
assumptions we notice that since both $f$ and $\Phi$ are norm
continuous it suffices to establish the inequality for the norm dense
set of $n-$tuples $\ux = (x_1, \dots , x_n)$, such that each $x_i$ has
finite spectrum. Since the $x_i$'s commute mutually there is in this
case a finite family $\{ p_k\}$ of pairwise orthogonal projections in
$A$ with sum $\bold 1$ such that 
$$
x_i = \sum_k \lambda_{ik} p_k\,, \qquad 1\le i \le n\,, 
\tag11
$$
where $\lambda_{ik} \in I_i$ (and repetitions may occur). If we set
$\ul_k = (\lambda_{1k}, \dots , \lambda_{nk})$ this means that
$$
f(\ux) = \sum f(\ul_k)p_k\,.  
\tag12
$$
As $\sum p_k = \bold 1$ in $A$ also $\sum \Phi (p_k)=1$ in $L^\infty_{\mu_C} 
(T)$, so that $\sum \Phi (p_k)(t) =1$ for (almost) every $t$ in
$T$. Since $f$ is convex this implies that
$$
\aligned
&f\bigl(\Phi (x_1)(t), \dots , \Phi (x_n)(t)\bigr) = 
f\left(\sum\lambda_{1k}\Phi (p_k)(t), \dots , \sum \lambda_{nk} \Phi (p_k)
(t)\right) \\
&\le \sum f\left(\lambda_{1k}, \dots , \lambda_{nk}\right) \Phi (p_k)(t) 
= \sum f(\ul_k)\Phi(p_k)(t)\,. 
\endaligned
\tag13
$$
Consequently, by (8), (13) and (6),
$$
\aligned
\varphi_C (f, \ux)\; &\le \int \sum f(\ul_k)\Phi (p_k)(t)\,d\mu_C (t) =
\sum f(\ul_k)\int \Phi (p_k)(t)\,d\mu_C (t)\\ 
&=\sum f(\ul_k)\tau (p_k) = \tau (\sum f(\ul_k)p_k) = \tau (f(\ux))\,.
\endaligned 
\tag14
$$

To prove the inequality for a general $C^*-$algebra $A$ with commuting
$C^*-$sub\-algebras $A_i$ consider the GNS representation $(\pi_\tau,
\frak H_\tau)$  associated with $\tau$. By construction,
cf\. \cite{{\bf 14}, 5.1.5},  we obtain a normal, semi-finite trace
$\widetilde\tau$ on the von Neumann algebra $M=\pi_\tau (A)''$ such
that  $\widetilde\tau( \pi_\tau (x))=\tau(x)$ for every $x$ in $A_+$,
and we can define  the two-sided ideal $M_\tau$ and its norm closure  
$M_\tau^=$ as in (10). If we now put  $\tA = M_\tau^=$ and 
$\tA_i = \tA \cap (\pi_\tau (A_i))^{-w}$ (closure in the weak operator
topology) for $1\le i\le n$, we have exactly the setup above. Our
argument, therefore, shows that the inequality holds 
in the setting of
$\tA,\, \tA_i$ and  $\widetilde\tau$, and an arbitrary commutative
$C^*-$subalgebra $\tC$ of $\tA$. 

Note now that since $\tau$ is densely defined on $A$ we have
$\pi_\tau(A)\subset \tA$. If  $C$ is a commutative $C^*-$subalgebra of
$A$ we put $\tC=\pi_\tau(C)$, and observe that $\tC$ has the form
$\tC=C_0(\tT)$ for some closed subset  $\tT$ of $T$ with 
$\mu_C (T\setminus \tT) =0$. The map $\widetilde \Phi$ from $\tA$ to
$L^\infty_{\mu_C} (\tT)$ defined in section 3 will therefore satisfy the 
restriction formula  
$$
\widetilde \Phi (\pi_\tau (x)) =  \Phi (x) \,|\, \tT
\tag15
$$
for every $x$ in $A$. But then, by our previous result in (14),
$$
\align
\varphi_C(f, \ux)\; & =\int f\left(\Phi (x_1), \dots ,
\Phi(x_n)\right)\,d\mu_C  
=\int f\left( \widetilde\Phi(\pi_\tau (x_1)), \dots , \widetilde\Phi
  (\pi_\tau (x_n)) \right)\,d\mu_C \\
&= \varphi_{\tC}(f, \pi_\tau (\ux)) \le \widetilde\tau (f(\pi_\tau(\ux)) 
=\widetilde \tau (\pi_\tau (f(\ux))= \tau (f(\ux))\,.
\tag16
\endalign
$$
\hfill$\square$
\enddemo

\bigskip 

\demo{Proof of Theorem 2} It is evident that  the function 
$\ux @>>> \varphi_C (f, \ux)$, defined in Lemma 6, is convex on 
$\bigoplus (A_i)^{I_i}_{\sa}$ for each $C$, being composed of a linear
operator  $\Phi (=\Phi_C)$, a convex function $f$, and a positive linear
functional  - the integral. Moreover, if $\Cal C$ denotes the set of
commutative $C^*-$sub\-algebras of $A$, it follows from Lemma 6 that  
$$
\tau(f(\ux)) = \sup_{\Cal C} \; \varphi_C (f, \ux )\,,
\tag17
$$
the supremum being attained at every $C$ that contains the commutative
$C^*-$sub\-algebra $C^*(\ux)$ generated by the (mutually commuting) elements 
$x_1, \dots , x_n$. Thus $\ux @>>> \tau(f(\ux))$ is convex as a supremum of 
convex functions. \hfill$\square$
\enddemo

This concludes the first part of our paper and we now turn to Theorem 3.

\bigskip 

\subhead{7. $\bold\ell-$Convexity}\endsubhead  We consider a strictly 
increasing, convex and continuous function $e$  on some interval
$I$, and denote by $\ell$ its inverse function (so that $\ell(e(s)))=s$ for
every $s$ in $I$ and $e(\ell(t))=t$ for every $t$ in $e(I)$) . If $g$ is
a convex function defined on a convex subset of a linear space then
$f=e\circ g$ is convex as well, but in some sense $f$ is ``much more''
convex, since even $\ell\circ f$ is convex, whereas $\ell$ is concave. We
say in this situation that the function $f$ is $\ell-${\it convex}. This
terminology is chosen to agree with the seminal example, where
$e=\exp$ and $\ell=\log$.  

Our aim is to show -- under mild restrictions on $\ell$ -- that
$\ell-$convexity has some remarkable structural properties, being preserved
under integrals and traces. By contrast, the concept of $\ell-${\it concavity}
-- with the obvious definition -- seems to be less interesting.

\bigskip

\proclaim{8. Lemma} If $I$ is an interval in $\Bbb R$ and $e$ is a
 strictly increasing, strictly convex function in $C^2(I)$ with
inverse function $\ell$, then for each probability measure $\mu$ on a locally
compact Hausdorff space $T$, and for each $\ell-$convex function $f$
defined on some cube $\uI$ in $\Bbb  R^n$, the function   
$$
(u_1, \dots ,u_n) @>>> \int f\left(u_1(t), \dots , u_n(t)\right)\,d\mu(t)\,, 
\qquad u_i\in L^\infty_\mu (T)\,,
\tag18
$$
is also $\ell-$convex on the appropriate $n-$tuples in $L^\infty_\mu (T)$ 
if and only if the function $\varphi$, defined by 
$\varphi(e(s))=(e''(s))^{-1}(e'(s))^2$, is concave.
\endproclaim

\demo{Proof}  Since $f=e\circ g$ for some convex function $g$ on
$\uI$, we see that to prove the Lemma it suffices to show that the increasing 
function  
$$
k(u)=\ell\left(\int e(u(t))\,d\mu(t)\right)\,, \qquad u\in L^\infty_\mu (T)\,,
\tag19
$$
is  convex on  $(L^\infty_\mu (T))^I$. Clearly, this is also a necessary
condition. Considering instead the scalar funtions
$$
h(s)=\ell\left(\int e\left(u(t)+sv(t)\right)\,d\mu(t)\right)\,
\tag20
$$
for arbitrary elements $u, v$ in $L^\infty_\mu (T)$, where the range
of $u$ is contained in the interior of $I$, we notice that convexity of $k$ is
equivalent to convexity at zero for all functions of the form $h$,  and we
therefore only have to show that $h''(0)\ge 0$. Setting  $r(s) =
\int e(u+sv)d\mu$ we compute 
$$ 
\aligned
&h'(s)= \ell'(r)\int e'(u+sv)v\,d\mu\,; \\
&h''(s)= \ell''(r)\left(\int e'(u+sv)v\,d\mu\right)^2
+\ell'(r)\int e''(u+sv)v^2\,d\mu\;.
\endaligned
\tag21
$$
On the other hand, since $\ell(e(t))=t$ we also have
$$
\aligned
&\ell'(e(t))e'(t)=1\,;\\
&\ell''(e(t))\left(e'(t)\right)^2 + \ell'(e(t))e''(t) =0\,.
\endaligned
\tag22
$$
In our case we can let $e(t) = r(s)$, whence $t =\ell(r(s))=h(s)$. We can
therefore eliminate $\ell''(r)$ in (21) to get the expression 
$$
\align
&h''(s)= -\ell'(r)e''(h(s))(e'(h(s)))^{-2}\left(\int e'(u+sv)v\,d\mu\right)^2
+ \ell'(r)\int e''(u+sv)v^2\,d\mu \\
&=\ell'(r)\left( \int e''(u+sv)v^2\,d\mu 
- e''(h(s))(e'(h(s)))^{-2}\left(\int e'(u+sv)v\,d\mu
\right)^2\right)\,. \,\quad \text{(23)}
\endalign
$$
Since $e$ is strictly increasing, so is $\ell$, which implies that
$\ell'(r)>0$. It follows from (23) that $h''(0)\ge 0$ if and only if 
$$
e''(h(0))e'(h(0))^{-2} \left(\int e'(u)v\,d\mu\right)^2 \le \int e''(u)v^2\,d\mu\,.
\tag24
$$
Now define the function $\varphi$ on $e(I)$ by 
$\varphi(e(s))=(e'(s))^2(e''(s))^{-1}$. Since $e(h(0))=r(0)=\int e(u)\,d\mu$ 
and $e''(h(0)) > 0$ because $e$ is strictly convex we see that (24)
is equivalent to the inequality
$$
\left(\int e'(u)v\,d\mu\right)^2 \le 
\varphi\left(\int e(u)\,d\mu\right) \int e''(u)v^2\,d\mu\,.
\tag25
$$
For ease of notation put $\varphi\left(\int e(u)\,d\mu\right)
=\widetilde\varphi(e(u))$, and choose a function $w$ in $L^1_\mu(T)$ with
$\int w\,d\mu =1$. Then consider the quadratic form
$$
\aligned
&\lambda^2 \int e''(u)v^2\,d\mu 
-2\lambda \int e'(u)v\,d\mu + \widetilde\varphi(e(u))\int w\,d\mu \\
=\; &\int\left( \lambda^2 v^2 e''(u) 
-2\lambda v e'(u) +\widetilde\varphi(e(u))w \right)\,d\mu\,.
\endaligned
\tag26
$$
By construction this form is positive if and only if (25) is satisfied. But
(26) expresses the integral of a function which is itself a quadratic
form. The minimum in (26) therefore occurs for $\lambda v e''(u)=
e'(u)$ and equals
$$
\aligned
&\int\left(\widetilde\varphi(e(u))w -
  (e''(u))^{-1}(e'(u))^2\right)\,d\mu \\
=\; &\int \left( \widetilde\varphi(e(u))w - \varphi(e(u))\right) \,d\mu =
\varphi\left( \int e(u)\,d\mu\right) -\int \varphi (e(u))\,d\mu\,.
\endaligned
\tag27
$$
Evidently this expression is non-negative if and only if $\varphi$
is concave.
\hfill$\square$ 
\enddemo

\bigskip 

\subhead{9. Remark}\endsubhead Using the equation $e(\ell(t))=t$ as in (22)
we easily find that $\varphi (t) = - (\ell''(t))^{-1} \ell'(t)$, so
that the condition in Lemma 8 translates to the demand:
$$ 
\text{The (negative) function}\quad t@>>> \ell'(t)/\ell''(t) \quad
\text{must be convex.}
\tag28
$$
Note also that the condition in (27), viz.
$$
\int \varphi(e(u))\,d\mu \le \varphi\left(\int e(u)\,d\mu\right)
\tag29
$$
makes sense for an arbitrary measure $\mu$ and provides a necessary
and sufficient condition for the function defined in (18) to be
$\ell-$convex. However, in order to satisfy (29) for an arbitrary (point)
measure, the function $\varphi$ must have $s\varphi (t)\le \varphi
(st)$ for all $s>0$, which forces it to be homogeneous, i\.e\. $\varphi
(st)=s\varphi(t)$ for $s>0$. It follows that the function defined by
(18) in Lemma 8 is $\ell-$convex for an arbitrary measure $\mu$ if and
only if 
$$
\ell'(t)/\ell''(t)= \gamma t \quad \text{for some non-zero number $\gamma$.}
\tag30
$$
Of course, this can only happen when the domain of $\ell$ is stable under
multiplication with positive numbers, so it is either a half-axis or the
full line. But since the expression in (30) must be negative, only
half-axes can occur.

\bigskip 

\demo{Proof of Theorem 3} The Theorem follows by using Lemma 6, as in the proof
of Theorem 2, combined with Lemma 8.
\hfill$\square$ 
\enddemo

\bigskip 

\subhead{10. Examples}\endsubhead Evidently the condition that $\ell'/\ell''$
be convex is not very restrictive, and is satisfied by myriads of functions,
of which we shall list a few, below. On the other hand, the demand that
$\ell'/\ell''$ be homogeneous is quite severe, and only four (classes of)
functions will meet this requirement:

\noindent{\bf (i)} Let $\ell (t) = \log (t)$ for $t>0$. We get $e(s)=\exp (t)$
for $s$ in $\Bbb R$ and $\ell'/\ell'' = -t$. This is the classical example,
and by far the most important. Clearly $\ell (t)= c\log (t)$ for any $c>0$ can
also be used, but we omit this trivial parameter here and in the following
examples.

\noindent{\bf (ii)} Let $\ell (t)= t^{1/p}$ for $t\ge 0$ and some $p>1$. We
get $e(s)=t^p$ for $s\ge 0$ and $\ell'/\ell'' = -\gamma t$, where $\gamma = 
p/(p-1)>1$. The root examples are also fairly well known. Indeed, it is a very
general fact that whenever $f$  is a convex (resp. concave) function that  is
homogeneous of some degree $p>0$ then (i) we must have $p\ge 1$ (resp. $p\le
1$) and (ii) the function $f^{1/p}$ is automatically  convex (resp.
concave). This is discussed in detail in the proof of Corollary 1.2 in
\cite{{\bf 9}}.

\noindent{\bf (iii)} Let $\ell (t)= - t^{-\alpha}$ for $t>0$ and some $\alpha
>0$. We get $e(s)= (-s)^{-1/\alpha}$ for $s <0$ and $\ell'/\ell'' = 
-\gamma t$, where $\gamma=(1+\alpha)^{-1} < 1$.

\noindent{\bf (iv)} Let $\ell (t)=-(-t)^p$ for $t\le 0$ and some $p>1$. We get
$e(s)= -(-s)^{1/p}$ for $s\le 0$ and $\ell'/\ell''=\gamma t$, where
$\gamma =(p-1)^{-1}>0$. 

Non-homogeneous examples are not hard to come by. Without any apparent order
we mention these:

\noindent{\bf (v)} Let $\ell (t)=-\exp(-\alpha t)$ for $t$ in $\Bbb R$ and
some $\alpha >0$. We get $e(s)=-\alpha^{-1}\log (-s)$ for $s <0$ and
$\ell'/\ell'' = -1/\alpha$.  In applications of Theorem 3 the parameter $\alpha$
disappears, since $\ell(\tau(e(a)))=-\exp(-\alpha\tau(-\alpha^{-1}\log(-a)))=
-\exp(\tau(\log(-a)))$ for any operator $a<0$; so we may as well assume that
$\alpha =1$.

\noindent{\bf (vi)} Let $\ell(t)=\log(\log(t))$ for $t>1$. We get 
$e(s)=\exp(\exp(s))$ for $s$ in $\Bbb R$ and $\ell'(t)/\ell''(t)=
-t\log(t)(1+\log(t))^{-1}$, which is only convex for $t\le e$. So on the 
intervals $1<t\le e$ and $-\infty <s\le 0$ we can use the functions 
$\log\log$ and $\exp\exp$.

\noindent{\bf (vii)} Let $\ell(t)=(\log(t))^{1/p}$ for $t\ge 1$ and some
$p>1$. Here $e(s)=\exp(s^p)$ for $s \ge 0$ and by computation we find that
$\ell'/\ell''$ is convex for $t\le \exp(1+1/p)$. The allowed interval for $e$
is $0\le s \le (1+p)p^{-2}$.

\noindent{\bf (viii)} Let $\ell(t)=\left(1-(1-t)^p \right)^{1/p}$ for 
$0\le t\le 1$ and some $p>1$. We get $e(s)=1-(1-s^p)^{1/p}$ for 
$0\le s \le 1$ and $\ell'(t)/\ell''(t)=(p-1)^{-1}\left((1-t)^{p+1}
-(1-t)\right)$, which is a convex function on the unit interval. 

\noindent{\bf (ix)} Let $\ell(t)=t^{1/p}(1+t^{1/p})^{-1}$ for $t\ge 0$ and
some $p\ge 1$. Here $e(s)=s^p(1-s)^{-p}$ for $0\le s <1$. For $p>1$ we find 
after some computation that $\ell'(t)/\ell''(t)=
2pt^{1+1/p}\left((p-1)(p-1+(p+1)t^{1/p})\right)^{-1}$ plus a linear term, 
and this is a convex function. For $p=1$ we simply get $\ell'(t)/\ell''(t)=
-\tfrac12 (1+t)$, and we note that this example is just a translation of (iii)
(replacing $t$ by $1+t$ and adding 1).

\bigskip 
  
\proclaim{11. Corollary} If $f$ is a positive, continuous, $\log-$convex
function on a cube $\uI$ in $\Bbb R^n$, and $A_1, \dots , A_n$ are 
mutually commuting  $C^*-$subalgebras of a $C^*-$algebra
$A$, then for each densely defined, lower semi-continuous trace $\tau$
on $A$ the function 
$$
(x_1, \dots , x_n) @>>> \tau (f(x_1, \dots , x_n))\,, 
\tag31
$$
defined on commuting $n-$tuples such that $x_i\in (A_i)_{sa}^{I_i}$ for each
$i$, is also $\log-$convex on $\bigoplus (A_i)_{\sa}^{I_i}$.
If, instead, $f^{1/p}$ is convex for some $p>1$, or if $-f^{-\alpha}$ is
convex for some $\alpha >0$, or if $f<0$ and $-(-f)^p$ is convex, then we also
have convexity of the respective functions  
$$
\aligned
&(x_1, \dots , x_n) @>>> (\tau (f(x_1, \dots , x_n)))^{1/p}\,,\\
&(x_1, \dots , x_n) @>>> -(\tau (f(x_1, \dots , x_n)))^{-\alpha}\,,\\
&(x_1, \dots , x_n) @>>> -(-\tau(f(x_1, \dots , x_n)))^p\,.
\endaligned
\tag32
$$
\endproclaim

\bigskip 

\subhead{12. Remarks}\endsubhead The Corollary above applies to some
unexpected situations. Thus we see that the function
$$
(x,\,y) @>>> \log(\tau(\exp((x+y)^2)))\,,
\tag33
$$
where $x$ and $y$ are self-adjoint elements in a pair of commuting
$C^*-$algebras $A$ and $B$, is (jointly) convex. The same can be said of the
function
$$
(x,\,y) @>>> \log(\tau(\exp(- x^\alpha y^\beta)))
\tag34
$$
for $0 < \alpha,\,\beta$ and $\alpha+\beta\le 1$, defined on $A_+\times B_+$. 
Applied to the root functions the Corollary shows that the function
$$
(x,\,y) @>>> \left(\tau((x+y)^q)\right)^{1/p}
\tag35
$$
is convex for $1\le p\le q$, and that also
$$
(x,\,y) @>>>  
\left(\tau\left((\bold 1 -x^{\alpha} y^{\beta})^p\right)\right)^{1/p}
\tag36
$$
is convex for $p\ge 1$ on the product of the positive unit balls of $A$ and
$B$. 

The last two cases in Corollary 11 are perhaps easier to apply in terms of
concave functions. By elementary substitutions we find that if $f$ is a
positive concave function on some cube $\uI$ in $\Bbb R^n$, then both
functions 
$$
\aligned
&\ux @>>> \left(\tau\left((f(\ux))^{-\alpha}\right)\right)^{-1/\alpha}\,;\\
&\ux @>>> \left(\tau\left((f(\ux))^{1/p}\right)\right)^p\,;
\endaligned
\tag37
$$
are concave on $\oplus (A_i)^{I_i}_{\sa}$ for $\alpha >0$ (and $f>0$), and for
$p\ge 1$. In particular we see that 
$$
\left(\tau((x+y)^{1/p}))\right)^p \ge \left(\tau(x^{1/p})\right)^p +
 \left(\tau(y^{1/p})\right)^p\,,
\tag38
$$
for all $x,\, y$ in $A_+$, so that the Schatten $p-$norms are
super-additive for $p<1$.

\bigskip

Recall from \cite{{\bf 4}} the definition of the {\it  Kadison - Fuglede 
determinant}  $\Delta$ associated with a tracial state $\tau$ on a
$C^*-$algebra $A$:
$$
\Delta (x) = \exp(\tau (\log |x|)) \quad\text{whenever}\quad x\in A^{-1}\,. 
$$
This is a positive, homogeneous and multiplicative map on the set of  
invertible elements, cf\. \cite{{\bf 4}, Theorem 1}, closely related to
the ordinary determinant from matrix theory. Note, however, that when
$A=\Bbb M_n (\Bbb C)$ then $\Delta (x) = \left(\det (|x|)\right)^{1/n}$, 
because we have to use the normalized trace $\tau = \tfrac1n\Tr$.

The following result is well known for matrices, at least for
functions of one variable. 
\medskip

\proclaim{13. Proposition} For each strictly positive concave function
$f$ on a cube $\uI$ in $\Bbb R^n$, and mutually commuting
$C^*-$subalgebras  $A_1, \dots , A_n$  of a $C^*-$algebra $A$ the
operator function  
$$
\ux @>>> \Delta (f(\ux))
\tag39
$$
is concave on the appropriate $n-$tuples of commuting self-adjoint
elements from  $\bigoplus (A_i)_{\sa}^{I_i}$. In particular, the
Kadison - Fuglede  determinant is a concave map on the set of positive
invertible elements. 
\endproclaim

\demo{Proof} From  example (v) in section 10 we see that the function
$$
\ux @>>> -\exp \left(\tau\left(\log (-(-f(\ux)))\right)\right) =
-\Delta (f(\ux))
\tag40
$$
is convex, as desired.   \hfill$\square$
\enddemo

\bigskip

\subhead{14. Remark}\endsubhead The concavity in (39) is closely related to
one of the main results (Theorem 6) in \cite{\bf 9} namely:
$x @>>> \tau \left(\exp (z +\log (x)\right)$ is concave for any 
self-adjoint operator $z$. In one sense (39) is stronger because it allows a
general concave $f$, but when $f$ is linear (39) is a corollary of Theorem 6
in \cite{\bf 9}, as we show now for one variable: 

We have to prove that \newline \centerline{$\exp\left(\tau \left(\log
(\tfrac12(x+y))\right)\right) \ge \tfrac1{2} \exp\left(\tau \left(\log
(x)\right)\right) + \tfrac1{2}\exp\left(\tau\left(\log
(y)\right)\right)$.}\newline
Put $z= -\log \left(\tfrac12(x+y)\right)$. Our condition is then that
$1 \ge \tfrac1{2} \exp\left(\tau \left(z+ \log (x)\right)\right) + \tfrac1{2} 
\exp\left(\tau \left(z+ \log (y)\right)\right)$. However, $\tau$ is a
state and therefore  Jensen's inequality applies. Thus,
$\exp\left(\tau \left(z+ \log (x)\right)\right) \le  \tau
\left(\exp\left(z+ \log (x)\right)\right)$, and similarly for $y$. By Theorem
6 of \cite{\bf 9} we know that $x @>>> \tau \left(\exp \left(z+ \log
(x)\right)\right)$ is concave, which implies that $\tfrac12\tau (\exp
(z+ \log (x))) + \tfrac12\tau (\exp (z+ \log (y)))  \le   \tau (\exp (z+ \log
(\tfrac12(x+y)))) = \tau (\bold 1) =1$, as desired.

\bigskip 

\subhead{15. Operator Means}\endsubhead We recall from 
\cite{{\bf 8}}, see also \cite{{\bf 6}, 4.1}, that a {\it Kubo-Ando mean}
on the set  $\bold B_+ =\Bbb B(\frak H)_+$ of positive operators is a
function  $\sigma  \colon  \bold B_+ \times \bold B_+ \to \bold B_+$
such that  

\noindent {\bf (i)} \quad $(\bold 1\;\sigma\; \bold 1) = \bold 1$, 

\noindent {\bf (ii)}  \quad $(x_1\;\sigma\; y_1) \le  (x_2 \;\sigma\; y_2)$ if
$x_1\le x_2$  and $y_1 \le y_2$,  

\noindent {\bf (iii)}  \quad $\sigma$ is (jointly) concave on $\bold B_+ \times
\bold B_+$, 

\noindent {\bf (iv)}  \quad  $z^*(x \;\sigma\; y) z =  ((z^*xz)\;\sigma\;
(z^*yz))$ for every invertible $z$ in $\Bbb B (\frak H)$. 

Of particular interest are the {\it harmonic mean} $!$ and the {\it
geometric mean} $\#$  defined on positive, invertible operators by:
$$
\aligned
&x\; ! \; y = 2(x^{-1}+y^{-1})^{-1}= 2x(x+y)^{-1}y\,;\\
&x\;\#\; y = x^{1/2}(x^{-1/2}yx^{-1/2})^{1/2}x^{1/2}\,.
\endaligned
\tag41
$$
($x\;\#\; y$ was introduced in \cite{{\bf 18}} and $\tfrac12 (x\; ! \; y)$,
the  {\it parallel sum}, was introduced in \cite{{\bf 1}}; see also \cite{{\bf
2}}, \cite {{\bf 6}, 4.1} and  \cite{{\bf 8}}.) The domain of
definition for these means can be  extended to all positive operators
by a simple limit argument (replacing $x$ and $y$  with
$x+\varepsilon\bold 1$ and $y+\varepsilon\bold 1$, and taking the norm
limit as $\varepsilon \to 0$), and we shall tacitly use this procedure
in the following. Thus we shall state all 
results for positive operators, but in the proofs assume that they are
invertible as well. Note that these two operator means are symmetric
in the two variables, and that they reduce to the classical harmonic and
geometric means when $x$ and $y$ are positive scalars, i\.e\.
$$
x\; ! \; y=2xy(x+y)^{-1} \quad \text{and} \quad x\;\#\; y =\sqrt{xy}.
\tag42
$$

A straightforward application of the Cauchy-Schwarz inequality shows that
$$
\tau (x\;\#\; y)\le \tau (x)\;\#\; \tau (y)
\tag43
$$
for any trace $\tau$. The corresponding result for the harmonic mean was
proved for the ordinary trace Tr in \cite{{\bf 2}}. The general version below
(Proposition 21) is somewhat more involved, but also richer. For greater
generality, but at little extra cost, we introduce the harmonic mean of an
$n-$tuple of positive operators as
$$
(x_1\;!\; x_2\; !\;\cdots \;!\;x_n)=n(x_1^{-1}+x_2^{-1}+\cdots
+x_n^{-1})^{-1}\,, 
\tag44
$$
and we note that this mean is symmetric in all the variables and increasing in
each variable. The fact that it is also jointly concave may not be widely 
known, so we present a short proof.

\bigskip

\proclaim{16. Proposition}  The $n-$fold harmonic mean is a jointly concave
function. 
\endproclaim

\demo{Proof} As the expression in (44) is homogeneous all we have to show is 
that
$$
\left(\sum x_i^{-1}\right)^{-1} + \left(\sum y_i^{-1}\right)^{-1}\le
\left(\sum (x_i + y_i)^{-1}\right)^{-1}\,,
\tag45
$$
for any pair of $n-$tuples of positive invertible operators. Multiplying left
and right by $\sum (x_i+y_i)^{-1}$ we obtain the equivalent inequality
$$
\aligned
\left( \sum (x_i+y_i)^{-1} \right) \left(\left(\sum x_i^{-1}\right)^{-1} +
\left(\sum y_i^{-1}\right)^{-1}\right)\left(\sum (x_i+y_i)^{-1}\right) \\
\le \sum (x_i + y_i)^{-1}\,.
\endaligned
\tag46
$$

We now appeal to the fact that the operator function $(x,\,y) @>>>
y^*x^{-1}y$ is jointly convex, \cite{{\bf 11}}, hence also jointly  
subadditive on the space of operators $x,\,y$, where $x$ is positive and
invertible. To see this, consider $n-$tuples $(x_i)$ and $(y_i)$ and define
$z_j = x_j^{-1/2}(y_j - a)$, with $a=\left(\sum x_i\right)^{-1}\sum y_i$. 
Then by computation we obtain the desired estimate 
$$
0\leq \sum z_j^* z_j = \sum y_j^* x_j^{-1}y_j -
\left(\sum y_j^*\right)\left(\sum x_j\right)^{-1}\left(\sum y_j\right)
\, .        
\tag47
$$

Breaking the left hand side of (46) into the sum of two terms and
using (47) on each we obtain the larger operator
$$
\aligned
&\sum (x_i+y_i)^{-1}x_i(x_i+y_i)^{-1}) + 
\sum (x_i+y_i)^{-1}y_i(x_i+y_i)^{-1}))\\ 
=\; &\sum (x_i+y_i)^{-1}(x_i+y_i)(x_i+y_i)^{-1} = \sum (x_i+y_i)^{-1}\,,
\endaligned
\tag48
$$
which is precisely the right hand side of (46), as claimed. \hfill$\square$
\enddemo

\bigskip 

\subhead{17. Remark}\endsubhead  Note that (47) also shows that the $n-$fold
harmonic mean is dominated by the arithmetic mean (the average).  Indeed,
$$
\align
& x_1\;!\;x_2\;!\;\dots \;!\;x_n 
= \left( \tfrac 1n \sum_{k=1}^n x_k^{-1}\right)^{-1}
=\left(\sum_{k=1}^n\tfrac 1n\bold 1\right)\left(\sum_{k=1}^n\tfrac 1n
x_k^{-1}\right)^{-1}\left(\sum_{k=1}^n\tfrac 1n\bold 1\right)\\  
\le \;&\sum_{k=1}^n\tfrac 1n\, \bold 1\,(x_k^{-1})^{-1}\,\bold 1
=\tfrac 1n\sum_{k=1}^nx_k\;.  
\tag 49
\endalign
$$
This result is not surprising, since the harmonic mean (on pairs of
operators) is the smallest symmetric mean, whereas the arithmetic mean
is the largest, cf\. \cite{{\bf 6}, 4.1}.

\bigskip

\proclaim{18. Proposition}  Given positive, invertible operators $x_1,x_2,
 \dots, x_n$ and $y$ in $\Bbb B (\frak H)$, let $d=diag\{x_1,x_2 \dots, x_n\}$
and $e=\sum_{i,j=1}^n y\otimes e_{ij}$ in $\Bbb M_n (\Bbb B (\frak H))$. Then
the following conditions are equivalent:
\roster
\item"{\bf (i)}"   \quad $y \le \left(\sum_{k=1}^n x_k^{-1}\right)^{-1}$\,,
\item"{\bf (ii)}"  \quad $ed^{-1}e \le e$\,,
\item"{\bf (iii)}" \quad $e \le d$\,.
\endroster
\endproclaim

\demo{Proof} Assume first that $y= \bold 1$, and set 
$a= \sum_{k=1}^n x_k^{-1}$, so that condition (i) becomes equivalent with  
$a\le \bold 1$. 

\noindent {\bf(i)}$\iff${\bf(ii)}. By computation
$$
ed^{-1}e=\sum_{i,j=1}^n a\otimes e_{ij},
\tag 50
$$
from which it follows that $a\le \bold 1$ if and only if $ed^{-1}e \le e$.

\noindent {\bf(i)}$\implies${\bf(iii)}. Let $p=\tfrac1n e$, and note that $p$
is a projection in $\Bbb M_n(\Bbb B(\frak H))$. For every $\varepsilon > 0$ we
have 
$$
\align
d^{-1} &= (p+(\bold 1-p))d^{-1}(p+(\bold 1-p))\le
(1+\varepsilon)pd^{-1}p+(1+\varepsilon^{-1})(\bold 1-p)d^{-1}(\bold 1-p)\\ 
&=(1+\varepsilon)n^{-2}ed^{-1}e+(1+\varepsilon^{-1})(\bold 1-p)d^{-1}(\bold 1-p)\,. 
\tag51
\endalign
$$
Now $a\le \bold 1$ by (1) and a fortiori $d^{-1}\le \bold 1$. Moreover, (i)$\implies$ (ii)
and thus we get
$$
d^{-1} \le (1+\varepsilon)n^{-2}e+(1+\varepsilon^{-1})(\bold 1-p) =
(1+\varepsilon)n^{-1}p+(1+\varepsilon^{-1})(\bold 1-p)\,.
\tag52
$$
Taking inverses this means that 
$$d \ge (1+\varepsilon)^{-1} np + \varepsilon (1+\varepsilon)^{-1}(\bold 1-p)\,,
\tag53
$$
from which the desired inequality follows as $\varepsilon \to 0$.

\noindent {\bf(iii)}$\implies${\bf(ii)}. If $e\le d$, then
$e+\varepsilon\bold 1 \le d+\varepsilon\bold 1$ for every $\varepsilon
> 0$, whence  
$$
(d+\varepsilon\bold 1)^{-1} \le (e+\varepsilon\bold 1)^{-1} = 
(n+\varepsilon)^{-1}p+\varepsilon^{-1}(\bold 1-p)\,.
\tag54
$$
But then, since $e(\bold 1-p)=0$, we have
$$
e(d+\varepsilon)^{-1}e \le (n+\varepsilon)^{-1}epe = (n+\varepsilon)^{-1}ne\,,
\tag55
$$
and as $\varepsilon \to 0$ we obtain the desired estimate.

In the general case, where $y$ is arbitrary (positive and invertible), we
define $\widetilde x_k = y^{-1/2}x_ky^{-1/2}$ for $1\le k \le n$. Then with
$\widetilde e = \sum_{i,j=1}^n \bold 1\otimes e_{ij}$ and $\widetilde 
d = diag\{ \widetilde x_1, \widetilde x_2, \dots, \widetilde x_n\}$ we 
observe that conditions (i), (ii) and (iii) become equivalent with the
conditions  
\roster
\item"{\bf (i) $\tilde{}$}" \quad $\bold 1\le \left(\sum_{k=1}^n 
\widetilde x_k^{-1}\right)^{-1}$\,,

\item"{\bf (ii) $\tilde{}$}" \quad 
$\widetilde e {\widetilde d}^{-1}\widetilde e\le \widetilde e $\,, 

\item"{\bf (iii) $\tilde{}$}" \quad $\widetilde e \le \widetilde d$\,.
\endroster
But these conditions are exactly the ones we proved to be equivalent above,
assuming that $y=\bold 1$. \hfill $\square$
\enddemo

\bigskip

\proclaim{19. Corollary}  For any $n-$tuple $(x_1, x_2, \dots , x_n)$ of 
positive operators in $\Bbb B(\frak H)$ the harmonic mean
$(x_1\;!\;x_2\;!\;\dots \;!\;x_n)$ is the largest positive operator of the
form $nz$, such that
$$
\pmatrix x_1 & 0 & \dots & 0\\0 & x_2 & \dots & 0\\  \vdots & \vdots &
\ddots & \vdots \\0 & 0 & \dots & x_n \endpmatrix 
\ge \pmatrix z & z & \dots & z\\z & z & \dots & z \\ \vdots & \vdots &
\ddots & \vdots  \\ z & z & \dots & z\endpmatrix
\tag56
$$
\hfill $\square$
\endproclaim

\bigskip

\subhead {20. Remark}\endsubhead  Note that Proposition 18 actually says
slightly more than Corollary 19, namely that the set of positive operators
$nz$, such that $z$ satisfies the matrix inequality in (56), is exactly equal
to
$$
\left((x_1\;!\;x_2\;!\;\dots\;!\;x_n)-\Bbb B(\frak H)_+\right) \cap 
\Bbb B(\frak H)_+\,.
\tag57
$$
This should be compared to the analogous result for the geometric mean of
pairs of positive operators $x$ and $y$. Here $x\;\#\;y$ is the largest
positive operator $z$ such that $\left( \smallmatrix x & z\\z &
y\endsmallmatrix\right) \ge 0$ in $\Bbb M_2(\Bbb B(\frak
H))$. However, one can find positive operators $z$ such that $z\le
x\;\#\;y$, without the matrix being positive. It suffices to take
$x=\bold 1$, so that $x\;\#\;y = y^{1/2}$. Evidently the matrix  
$\left(\smallmatrix \bold 1 & z\\z & y\endsmallmatrix\right)$ is positive only
if $z^2 \le y$, and it is easy to find examples where $z \le y^{1/2}$, but
$z^2\not\le y$. 

\bigskip

\proclaim{21. Proposition} If $\tau$ is a densely defined, lower
semi-continuous trace on a $C^*-$algebra $A$, and $f$ is a strictly positive
function on $]0,\infty[$ such that $t@>>> \left(f(t^{-1})\right)^{-1}$ is  
concave, then for all positive operators $x_1\,\dots , x_n \in A$  and 
$\alpha > 0$ we have 
$$
\aligned
&\left(\tau \left((f(x_1\; ! \; x_2\; !\;\cdots\; ! \;
x_n))^\alpha\right)\right)^{1/\alpha} \\
\le\; &\left(\tau((f(x_1))^{\alpha})\right)^{1/ \alpha}\; ! \;
\left(\tau(f((x_2))^\alpha)\right)^{1/\alpha} \; ! \; \cdots \; !\;
\left(\tau((f(x_n))^\alpha)\right)^{1/\alpha}\,.  
\endaligned
\tag58
$$
\endproclaim

\demo{Proof} Define $g(t)=(f(t^{-1}))^{-1}$. Then from (37) we see that
the function  
$$
x@>>> \left(\tau \left((g(x))^{-\alpha}\right)\right)^{-1/\alpha}
\tag59
$$
is concave on the set of positive invertible elements in $A$. In particular, 
$$
\aligned
&\left(\tau\left((g(1/n\, (x_1^{-1}+ \cdots
+x_n^{-1})))^{-\alpha}\right)\right)^{-1/\alpha}\\ 
\ge\; &1/n \left(\tau \left((g(x_1^{-1}))^{-\alpha}\right)\right)^{-1/\alpha}
+\cdots +  
1/n\left(\tau \left((g(x_n^{-1}))^{-\alpha}\right)\right)^{-1/\alpha}\,.
\endaligned
\tag60
$$
Taking the reciprocal values and noting that
$1/n \, (x_1^{-1}+\cdots + x_n^{-1}) = (x_1\; ! \; \cdots \; ! \; x_n)^{-1}$ 
this means that
$$
\aligned
&\left(\tau\left((g((x_1\; ! \; \cdots \; ! \;
x_n)^{-1}))^{-\alpha}\right)\right)^{1/\alpha}\\ 
\le \; &n\left(\left(\tau ((g(x_1^{-1}))^{-\alpha})\right)^{-1/\alpha}
+\cdots +\left(\tau ((g(x_n^{-1}))^{-\alpha})\right)^{-1/\alpha}\right)^{-1}\\
=\; &\left(\tau((g(x_1^{-1}))^{-\alpha})\right)^{1/\alpha}\;  ! \; \cdots
\; ! \;  \left(\tau ((g(x_n^{-1}))^{-\alpha})\right)^{1/\alpha}\,. 
\endaligned
\tag61
$$
Since $g(t^{-1})^{-1}=f(t)$ this is the desired result. \hfill $\square$
\enddemo

\bigskip 

\subhead{22. Remark}\endsubhead The result above applies to the
functions $t @>>> t^{1/p}$ for $p\ge 1$, in particular to the identity
function (and with $\alpha = 1$); but it also applies to the
functions $t @>>> \log (1+t^{1/p})$ for $p \ge 1$. 
 
On the abstract level Proposition 21 applies to
$C^2-$functions $f$, such that $f$ and $t @>>> t^2f(t)^{-2}f'(t)$ are
simultaneously increasing or decreasing (strictly) on $]0, \infty[$.

Proposition 21 also applies  to other (not necessarily  symmetric)
means of positive operators. Such means were introduced in 
\cite{{\bf 8}} (see also \cite{{\bf 6}, 4.1}).  

\bigskip

\proclaim{23. Proposition}  For every Kubo - Ando mean $\sigma$ and for
every densely defined, lower semi-continuous trace $\tau$ on a
$C^*-$algebra  $A$ we have $\tau (x \;\sigma\; y) \le \tau (x) \;\sigma\; \tau
(y)$ for all $x, y$ in $A_+$. 
\endproclaim

\demo{Proof}  For each mean $\sigma$ there is a unique probability
measure $\mu$ on $[0, \infty]$ such that  
$$
x\; \sigma\; y =\tfrac12 \int_0^\infty (tx\; !\; y)\,  (1+1/t)\,d\mu(t)
=\tfrac12 \int_0^\infty \left((1+t)x\; !\; (1+1/t)y\right)\,d\mu(t)\,,
\tag62
$$
cf\. \cite{{\bf 6}, 4.1}. Applying Proposition 22 with $f(t)=t$ and
$\alpha=1$ it follows that 
$$
\aligned
&\tau (x\; \sigma\; y) 
= \tfrac12 \int_0^\infty \tau (tx\; ! \; y)\,  (1+1/t)\,d\mu(t)\\
\le\; &\tfrac12 \int_0^\infty \left(t\tau (x)\; ! \; \tau(y)\right)\,(1+1/t)
\,d\mu(t) = \tau (x)\; \sigma\; \tau (y)\,.
\endaligned
\tag63
$$
\hfill$\square$ 
\enddemo

\bigskip

{\it Acknowledgements:} The work of E.H.L. was partly supported by 
U.S. National Science Foundation grant PHY98-20650-A02, whereas
G.K.P. was partly supported  by the Danish Research Council (SNF).  

\newpage

\enddocument